# NECESSITY AND CHANCE:
## DETERMINISTIC CHAOS IN ECOLOGY AND EVOLUTION

### ROBERT M. MAY

Abstract. This is an outline of my Gibbs Lecture to the American Mathematical Society in January 1994; it is essentially a sign-posted guide to a still-developing literature.

## Introduction

During the past twenty years or so, "deterministic chaos" has moved centre stage in many areas of applied mathematics. One important stimulus for this, particularly in the early 1970s, was work on nonlinear aspects of the dynamics of plant and animal populations. Here, I first sketch some of this biological background and then move on to recent work on chaos-based ideas for nonlinear forecasting in fluctuating time-series. Measles data from England and Wales as a whole, and from individual large cities, will be used to illustrate some of these ideas, especially those pertaining to relations between the spatial scale over which data are aggregated and our ability to perceive underlying patterns produced by largely deterministic mechanisms.

I then go on to show how very simple and fully deterministic mathematical models which incorporate rules for local movement on some spatial lattice can generate an extremely diverse array of spatial patterns, including spiral waves, apparently static but inhomogeneous "crystal lattices", and spatial chaos. I will outline examples ranging from the population dynamics of host-parasitoid and other systems, to evolutionary games (Prisoner's Dilemma, Hawk-Dove) played among "territory holders", and will speculate on the implications of this work and on its likely future directions.

In all this, I will be offering a sign-posed guide to the literature and not a treatment of any topic in depth. My aim is to persuade the mathematically oriented readers of this journal that mathematical biology is a discipline which has come of age. No longer dominated by routine applications of differential equations (or, at best, computer simulations of partial differential equations), the subject has—most notably in the case of chaos [1]—helped to advance new areas of mathematics. The subject is increasingly characterised by substantial engagement between mathematical models and data, often with practical implications [2, 3]. Even in relatively conventional areas of mathematics, recent work on the evolution of virulence [4], or on multiple-epitope models of the interactions between viruses and the immune system [5], or on the effects of habitat destruction in spatially distributed multi-species communities [6] uses a combination of analytic techniques and numerical studies to show how nonlinearities can produce extraordinarily rich and surprising dynamics, complete with threshold effects and self-organised complexity on many different time scales.

## Deterministic chaos

In the 1950s and 1960s, several ecologists studied first-order difference equations as deliberately oversimplified metaphors for the effects of density-dependent regulation upon single-species populations. In particular, Moran [7] studied the equation

$$(1) \qquad N_{t+1} = R_0 N_t \exp(-aN_t) \,,$$

in an entomological context, and Ricker [8] in a fisheries one. Both were captives of the mind-set of their time, looking for stable points and finding them. Both also found cycles, and suggestions of some form of irregular behaviour, in numerical studies. Moran (personal communication) even set a graduate student to study this odd behaviour, but the student never completed his degree and the question was not pursued. At the same time, as surveyed fairly fully in my 1976 review [1] (which, however, overlooked Sharkovsky [9]), various mathematicians—of whom the Finn Myrberg seems the first—were aware of the complicated behaviour exhibited by first-order difference equations with a hump (negative Schwarzian derivative) of tunable severity but failed to recognise the wider implications of such phenomena.

It was left to Li and Yorke [10, which coined the word "chaos"], Oster and me [11], and others independently to discover these features of first-order nonlinear difference equations. As Yorke has put it, we were not first, but by fully appreciating the broader implications of what is now called chaotic behaviour, we were "the last to discover these facts." The canonical example of such an equation is the now-familiar quadratic map,

$$(2) \qquad x_{t+1} = ax_t(1 - x_t).$$

Its properties can be fairly easily laid bare, using techniques more readily accessible than those needed, for example, to develop the differential calculus. As I have urged elsewhere [1], students' intuition would be enriched if these equations were met very early in their mathematical education, which I think is still too much governed by accidents in the development of mathematical history. Then, later on, students could encounter the much more complicated examples of chaos than arise in differential equations, where at least three dimensions are needed [12].

At first sight, the chaotic trajectories of eq. (2) with $4 \geq a > 3.57\ldots$, or of eq. (1) with $R_0 > 14.767\ldots$, look like sample functions of some random process. Second thoughts might suggest that, apparent unpredictability notwithstanding, we can indeed make confident projections in such chaotic systems if we know the underlying simple rule; furthermore, such rules can be reconstructed by plotting $x_{t+1}$ against $x_t$ for these 1-dimensional systems. If so, then such deterministic chaos is more apparent than real.

But a third look reveals the extreme sensitivity to initial conditions in maps such as those of eqs. (1) and (2). This sensitivity to initial conditions, which leads to effective unpredictability even if we know or can infer the rule, is indeed the defining characteristic of "deterministic chaos". The phenomenon can be made more precise by defining a "Lyapunov exponent", $\lambda$, for the map $x_{t+1} = F(x_t)$ :

$$(3) \qquad \lambda \equiv \langle \ln |dF/dx|_{x_t}\rangle.$$

Here, $x_i$ are iterates of the map, with the magnitude of the slope of the map being evaluated at each point on the trajectory; $\lambda$ then represents the average value of the natural logarithm of this slope, evaluated over a long run of iterations. The reason for this definition becomes clear if we consider the evolution of two trajectories whose initial points, $x_0$ and $x_0 + \varepsilon$, differ only by a small amount, $\varepsilon << 1$. The first iterates of the two trajectories will differ by

$$(4) \qquad F(x_0 + \varepsilon) - F(x_0) \simeq \varepsilon (dF/dx)_{x_0} + \mathcal{O}\varepsilon^2.$$

If we now take the crude approach of replacing $(dF/dx)$ by its long-term average value of $\exp(\lambda)$, we observe that after $T$ iterates the two trajectories, initially differing by $\varepsilon$, have diverged by $\varepsilon \exp(\lambda T)$. Thus, if $\lambda > 0$, initially adjacent trajectories will, on average, diverge further and further. After a time of order

$$(5) \qquad T_\lambda \simeq \ln(1/\varepsilon)/\lambda,$$

the trajectories will be utterly divergent and all predictability will be lost. Since the error or uncertainty in initial conditions $\varepsilon$ only enters logarithmically, we have the order-of-magnitude estimate $T_\lambda \sim 1/\lambda$. We could call $T_\lambda$ the "Lyapunov horizon" for predictions in such deterministically chaotic systems: roughly speaking, for $T < T_\lambda$, knowledge of the underlying deterministic rule and the initial conditions will enable a degree of predictability, diminishing as $T$ increases toward $T_\lambda$; beyond the Lyapunov horizon, sensitivity to initial conditions renders prediction impossible, even though we are dealing with a fully deterministic system and we know the rules! For a more detailed but still relatively accessible account, see Mullin [13].

Pure mathematicians may be somewhat pained by the rather cavalier treatment of averages and other points of rigor. As will be seen in the next section, such concern can be justified, and much recent work focuses on questions about "local Lyapunov exponents" rather than broad-brush averages [14, 15]. But I make no apologies for the intuitive presentation above, which is correct in its essentials. Too much work in this area is presented in a way that is formal to the point of obscurantism.

## Chaos and forecasting

Given that very simple, and fully deterministic, rules or equations can produce time series that superficially look like random noise, we must look at all such apparently noisy time series with fresh eyes. This is, as it were, the flip side of the chaos coin. Applied to fluctuating and irregular time series such as numbers of voles [16], or cases of measles [17, 18], or cardiac rhythms [19], or economic time series of various kinds [20], the question is whether we are looking at "real randomness" (in the sense of a roulette wheel or other high-dimensional complexity) or a low-dimensional but chaotic attractor.

If we have a very long run of data, along with the assurance that there are no secular trends in the parameters of the underlying dynamical system, then various scaling-law techniques can be used to assess the dimensionality of the system and thence to identify low-dimensional deterministic chaos. For a good account of the essential ideas, and of some pitfalls, see Ruelle [21].

For most examples of practical interest, however, we have at least one of three kinds of problems: (a) the time series is not long enough for the data-hungry scaling law techniques to be applied; (b) there may be secular changes in the systems' parameters (we may indeed be dealing with the chaotic regime of a forced pendulum, but the pendulum length keeps changing); (c) the low-dimensional chaotic dynamics may be combined with sampling error, or "real randomness" (in the above-defined sense). The result is a rapidly expanding body of work which, to varying degrees, combines insights with frank phenomenology [21–25].

Here, I briefly sketch one approach to these problems and illustrate some of the surprises that can emerge. The essential ideal in this approach is to make short-term predictions that are based on a library of past patterns in a time series. By comparing the predicted and actual trajectories, we can then make tentative distinctions between dynamical chaos versus white noise or measurement error: for a low-dimensional chaotic time series, the accuracy of the nonlinear forecast falls off with increasing prediction-time interval (at a rate which gives an estimate of the Lyapunov exponent), whereas for uncorrelated noise, the forecasting accuracy is roughly independent of prediction interval. For a relatively short time series, it is more difficult to distinguish between low-dimensional chaos and autocorrelated or "coloured" noise; despite some tentative suggestions, there are unresolved problems here, not least because some forms of high-dimensional chaos may be equivalently viewed as low-dimensional chaos plus coloured noise [24, 26].

As developed by Farmer and Sidorowich [27], Sugihara and May [26] and later workers, this chaos-motivated approach to nonlinear forecasting first constructs a library of past patterns. This involves choosing an "embedding dimension", $E$, which uses time-lagged coordinates to represent each lagged sequence of data points $\{x_t, x_{t-\tau}, \ldots, x_{t-(E-1)\tau}\}$ as a point in an $E$-dimensional space. For each sequence for which we wish to make a prediction—each "predictee"—we now find those $E+1$ points in the library of $E$-dimensional points which are closest to the predictee in some sense (e.g., the vertices of the surrounding simplex which has least volume, or the minimum diameter). The prediction is now obtained by projecting the domain of the simplex into its range; that is, by keeping track of where the vertices of the simplex end up, $p$ time steps into the future. To get the prediction, we compute where the original predictee has moved within the range of the simplex, giving exponential weight (or something similar) to its original distances from the relevant neighbours. Finally, we compare the resulting prediction with the observed value, for each of a large number of such predictions, $p$ time steps into the future. To quantify the average reliability of the prediction thus arrived at, we compute the conventional statistical coefficient of correlation $\rho$ (or some equivalent measure) between predicted and observed values, as a function of the prediction time interval, $T_p$ (or the number of time steps $p$); $\rho$ is then plotted as a function of $T_p$.

There are clearly many empirical elements in this process. First, we have to guess the embedding dimension, $E$, and the appropriate lag, $\tau$, that combine effectively to relate our 1-dimensional time series to the higher-dimensional system that underlies it (for a more complete discussion, see [21, 28]). In the phenomenological approach described above, $E$ and $\tau$ are regarded as adjustable parameters, and we choose those values that give the highest values of $\rho$. This contrasts with Grassberger-Procaccia [29] and other scaling-law approaches, which seek to infer the embedding dimension directly. The tricky problem of secular trends in the underlying parameters (which is always a possibility for real data, in contrast with the "toy models" that constitute so large a fraction of the formal literature of this subject) can be brought into the picture by asking how we choose our "out of sample" predictions. For the measles data illustrated below, Sugihara et al. [17] have variously used the first half of the time series to compile the library and the second half to make out-of-sample predictions, or vice versa, or the entire time series to compile the library and then out-of-sample predictions for all time points, but with five years around each predictee removed from the library. The results for the measles data are insensitive to which procedure is used, encouraging the belief that secular trends in parameters are not a problem in this context.

Figure 1 shows the fit between predicted and observed cases of measles in England and Wales, as measured by the correlation coefficient $\rho$ as a function of prediction interval $T_p$ in monthly time steps. Here the nonlinear forecasting method does no better than conventional piecewise linear autoregressive techniques. The pattern seems to be primarily one of annual periodicity together with a roughly 2-year "interepidemic cycle", with the remaining variance being attributable to noise. This is surprising, because earlier analyses of similar prevaccination data for measles in New York City showed a marked decline in $\rho$ with increasing $T_p$, suggesting low-dimensional chaos [18, 26]. The answer to this puzzle lies in the degree of aggregation of the data. Figure 2 is as for Figure 1, except now the data are measles cases in Birmingham rather than England and Wales overall. In direct contrast with Figure 1, Birmingham shows the predictability declining with increasing prediction interval, characteristic of chaos as the Lyapunov horizon is approached.

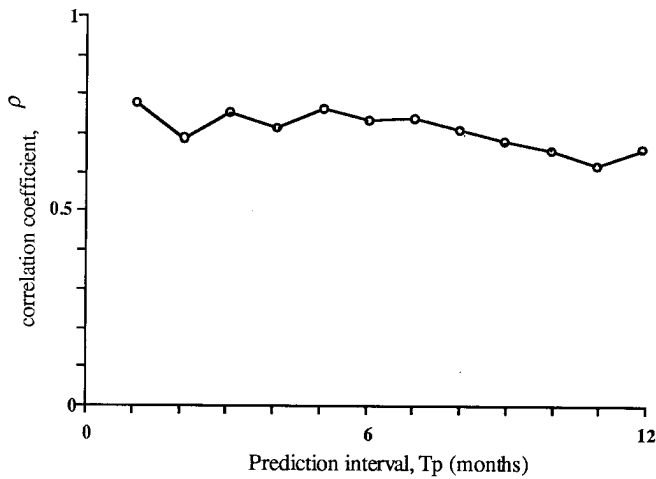

FIGURE 1. Predictability (as measured by the coefficient of correlation $\rho$ between predicted and observed numbers of measles cases) as a function of prediction interval (number of months ahead, $T_p$), for measles in England and Wales, 1948–1966. Here, the first half of the time series was used to construct the "library" of past patterns, which then was the basis for nonlinear forecasting for the second half of the data. The optimal embedding dimension and spacing/lag were chosen phenomenologically: $E = 8$, $\tau = 2$. The relatively flat $\rho$ versus $T_p$ relation suggests rough periodicity plus additive noise, and not low-dimensional chaos (after [17]).

There are further refinements to this measles story [30], but the main features are clear and instructive. Basic biological understanding of the population biology of the interaction between measles and its human hosts in endemic situations suggests we might expect the dynamics of the forced pendulum, with the "pendulum period" deriving from the host-parasite dynamics (corresponding to an intrinsic tendency to oscillation with a roughly 2-year period) and the forcing coming from the annual cycles of the school year [3]; such a system can exhibit complicated, almost periodic, chaotic dynamics. When we look at epidemiologically homogeneous units, in the form of cities big enough to maintain measles endemically, nonlinear forecasting techniques indeed suggest chaotic dynamics (and also give short-term predictions which are much better than the best linear predictions). But if we aggregate the data too coarsely, by looking at all measles cases in England and Wales, we lose the signature of the chaotic pendulum and see only periodicity plus noise (with nonlinear predictors doing neither better nor worse than linear ones; the individual, and short-term predictable, kinks and wiggles that ensue from the chaotic dynamics in any one city are all averaged out in the aggregated data set). All this highlights the extent to which answers about dynamical mechanisms depend on examining questions at the right spatial and temporal scale.

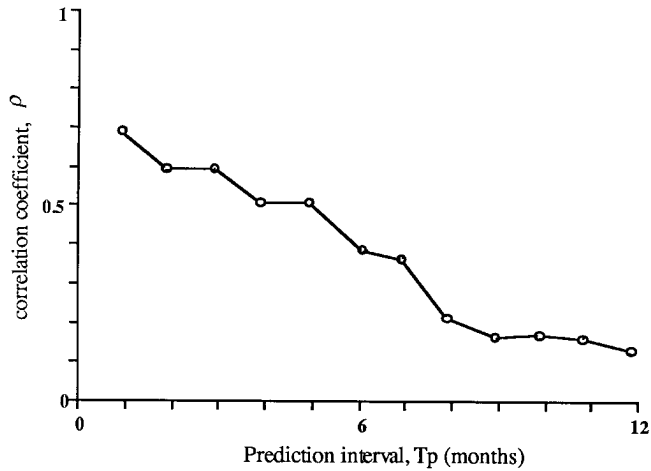

FIGURE 2. As for Figure 1, except now the data are for measles cases in one large city, Birmingham ($E = 7$, $\tau = 2$). In contrast with Figure 1, the $\rho$ versus $T_p$ relation shows a decline in predictability with increasing prediction interval characteristic of low-dimensional chaos. Similar results are found for other large cities in the UK and in the USA; London gives results intermediate between Figures 1 and 2 [17].

I have chosen here to focus on an example with real data, thus facing the messiness of the natural world. I could instead have used the more precise and elegant example of the quadratic map of eq. (2). This map, with $a = 4$, was used by Ulam and Von Neuman to generate random numbers on the MANIAC computer at the Institute for Advanced Study in Princeton in 1948 (no foreshadowing of chaos was involved in this choice). In this limit, eq. (2) gives numbers that appear randomly distributed on $[0, 1]$, with a probability density proportional to $[x(1-x)]^{-1/2}$. These pseudorandom numbers, suitably transformed, indeed are indistinguishable from white noise by conventional autoregressive linear tests (such as Bartlett's Kolmogorov-Smirnov test [26]). But the nonlinear forecasting methods described above clearly identify such an apparently random time series as low-dimensional chaos and indeed can predict the next one or two time points with $\rho = 0.89$. As we combine more and more coupled quadratic maps, however, the nonlinear forecasting technique cannot sort out the ensuing high-dimensional dynamics, and the system indeed appears random [17].

If a time series is generated by low-dimensional chaos, are some parts of the series more predictable than others? This question is of obvious interest, especially in applications to economic time series, such as currency exchange rates. Going back to eq. (4), the question can be rephrased in terms of *local* Lyapunov exponents. There may be some regions of the chaotic attractor for which adjacent points diverge significantly more slowly than for other regions, in which case the globally averaged $\lambda$ of eq. (3) will be less useful than a more fine-grained approach. Olsen and Schaffer [18; see their Figure 3] suggest evidence for just such regions of greater predictability in their analysis of measles data. For recent reviews of these important ideas, see Wolff [15] and Smith [14].

## Spatial chaos and population dynamics

Essentially all existing work on chaos deals with temporal aspects of the subject, either in a context where spatial dimensions are irrelevant or assuming spatial homogeneity. As outlined below, however, simple and fully deterministic rules for local movement, for individuals or groups within an overall population in an intrinsically homogeneous environment, can lead to self-organised spatial patterns, including spiral waves, apparently static but inhomogeneous "crystal lattices", and spatial chaos. Earlier work on spatially inhomogeneous systems in mathematical biology has largely dealt with reaction-diffusion equations (finding static patterns or travelling waves [31]) or with movement patterns in which individuals or groups were equally likely to move anywhere ("island" models, leading to relatively straightforward spatial dynamics [32]).

Here I discuss two explicit examples of simple, homogeneous, and deterministic systems which exhibit rich spatial dynamics, including spatial chaos. The first derives from a particular kind of prey-predator system of interest to ecologists. The second, in the next section, generalises the classic Prisoner's Dilemma and other game theoretic metaphors for evolutionary problems to include spatial considerations; the results are both mathematically interesting and aesthetically pleasing.

Parasitoids are insects (dipterans or hymenopterans) which lay their eggs in or on the larvae or pupae of arthropod hosts (usually lepidopterans). The emerging parasitoid offspring kill the host. These life history characteristics imply that the dynamical engagement between hosts and parasitoids is simpler than most other host-parasite or prey-predator associations (where the relation between parasitism of, or predation upon, the host and the subsequent production of parasite or predator progeny is much less direct). Mathematical models for host-parasitoid interactions can thus be brought into closer contact with field and laboratory data than is commonly the case for exploiter-victim systems, which helps explain the attention such systems have received [32]. Nor are parasitoids arcane or uncommon creatures; they probably constitute around 10% of all metazoan species, and they include many natural enemies of crop and orchard pests.

Specifically, let $H_t$ and $P_t$ represent the number (or density) of hosts and parasites in generation $t$, respectively. Then we may summarise the biological statements above by writing

$$(6) \qquad\qquad H_{t+1} = R_0 H_t F(P_t, H_t),$$

$$(7) \qquad\qquad P_{t+1} = cH_t[1 - F(P_t, H_t)].$$

Here $R_0$ is the intrinsic reproductive rate of the hosts (the number of offspring produced, on average, by each host in the absence of parasitoids—which may in general be density-dependent), $c$ is the number of (female) parasitoids to emerge from a parasitised host, and $F$ represents the probability that a host individual escapes parasitism (which, in general, will depend on the density of host and parasitoid populations in that generation). All the biological complications, which obviously can be substantial, are thus subsumed in the "searching function", $F(P, H)$.

The simplest assumption is that parasitoids search independently and randomly, in a spatially homogeneous world. The searching function is then the zeroth term in a Poisson distribution: $F(H, P) = \exp(-aP)$, where $a$ measures the parasitoid searching efficiency or attack rate. As shown by Nicholson and Bailey [33], the resulting dynamics are ever-diverging oscillations, resulting in extinction first of the hosts and consequently of the parasitoids.

So, what factors contribute to the persistence of host-parasitoid associations in the natural world? In laboratory studies, it seems that "mutual interference" among searching parasitoids may be important [32, 34]. But in field situations, the emerging consensus is that spatial patchiness, combined with a sufficient degree of variation in levels of parasitism among patches, may be the essential mechanism that enables the persistence of host and parasitoid populations at roughly steady overall densities [35, 36].

Given that sufficient variation in parasitoid attack rate among patches will result in overall persistence of the host-parasitoid association, the next question is what causes such variation? If the mechanism is an attack rate that depends directly or inversely on host density, then we must ask what produces the variability in host density from patch to patch. And if the mechanism is variability in attack rates uncorrelated with host density (possibly associated with variability in parasitoid densities among patches), then we similarly ask what causes such variability. In all cases, the usual answer is that external environmental factors cause patches to differ, thus producing variability in host and/or parasitoid densities among patches. An alternative possibility, however, is that such spatial structure may be self-organised, arising from simple and purely deterministic rules about the local movement of hosts and parasitoids within an environmentally homogeneous world, in which all patches are intrinsically identical [37–39].

Specifically, consider a 2-dimensional cartesian (chess board) array of $n \times n$ patches. In each patch, in each generation, hosts and parasitoids interact according to Nicholson-Bailey dynamics; that is, eqs. (6) and (7) with the random search function $F = \exp(-aP)$. This produces the next generation of hosts and parasitoids in that patch. Then, before the next round of searching and parasitism in each individual patch, a specified and density-independent fraction of the hosts ($\mu_H$) and of the parasitoids ($\mu_P$) move to the eight adjacent patches, being evenly distributed among them. In each patch, hosts and parasitoids now again interact, and so on.

This simple and biologically sensible model, with the dynamics in each patch being intrinsically unstable and with deterministic patterns of local movement, exhibits a truly remarkable range of dynamical behaviour. The densities of the host and parasitoid subpopulations in the 2-dimensional array of patches may exhibit complex patterns of spiral waves or spatially chaotic variation, or they may show static "crystal lattice" patterns, or they may become extinct. Figure 3 on the next page illustrates this range of dynamical behaviour. To do full justice to the dynamics, however, one really needs a video.

The various kinds of behaviour may be roughly understood. If the wave length of the spatio-temporal oscillations in host-parasitoid densities—the wave length of the spiral wave seen clearly in Figure 3a—is long compared with the side-length, $n$, of the arena, then we find the familiar Nicholson-Bailey extinction via ever-more-severe oscillations. Once the wave length (measured in relation to the numbers of patches) is smaller than $n$, so that a spiral wave can be "fitted into" the $n \times n$ arena, then hosts and parasitoids persist in the spiral wave patterns of Figure 3a. As the wave lengths get smaller and smaller in relation to $n$, we get more and more spiral waves within the arena; eventually there are sufficiently many that their patterns of collision and coalescing give rise to an ever-shifting jumble—spatial chaos. And in one corner of parameter space (with low host movement, high parasitoid movement, and $R_0$ not too big) we find static "crystal lattice" patterns, as illustrated in Figure 3c. Analytic results support the finding that these "crystal lattice" patterns need not be neatly periodic in space; they can be deformed, spatially irregular lattices, while retaining their essential property of being temporally static (these static lattices are only possible if $R_0 < e$).

I emphasise that this range of behaviour is obtained with the local dynamics being deterministically unstable, with a constant host reproductive rate, and with no density dependence in the movement patterns. The nature of the dynamics depends only on the host reproductive rate, $R_0$, on the values of the movement parameters, $\mu_H$ and $\mu_P$, and, to a degree, on the size of the arena, $n$ [37, 38].

These results, moreover, are relatively insensitive to the details of the interactions (i.e., the form of $F$) within any one patch. Thus we get essentially the same results from mathematically explicit models for the host-parasitoid interactions, such as Nicholson-Bailey (Figures 3a–c) or so-called Lotka-Volterra interactions (Figure 1d), and from very general "cellular automaton" models in which only qualitative rules are specified. Solé and Valls [39] have also obtained similar results, using corresponding models in continuous—as distinct from our discrete—space and time; they find results essentially identical to those in Figure 3.

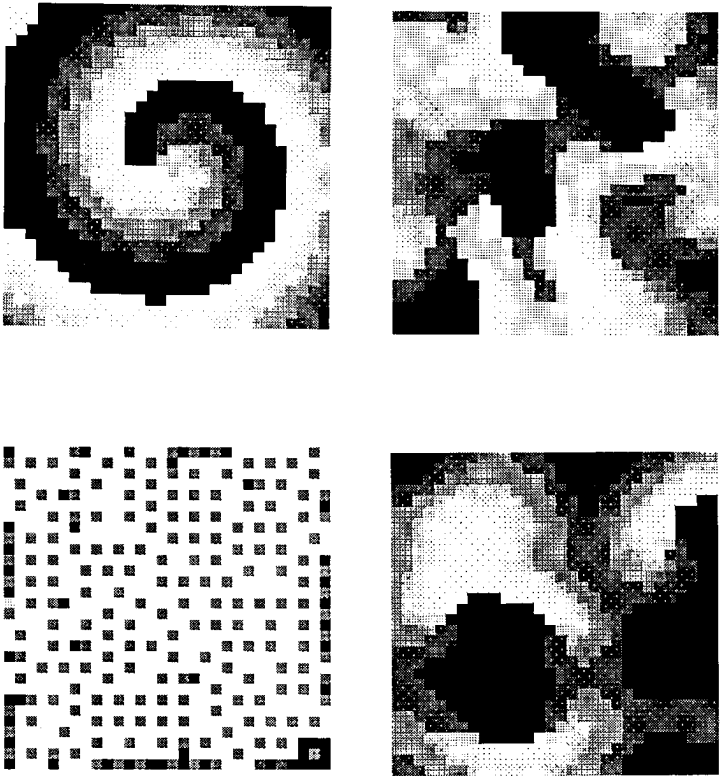

Figure 3. Instantaneous maps of the population density of hosts and parasitoids for simulations based on the model described in the text, with local dispersion and Nicholson-Bailey local dynamics. Here $\lambda = 2$ and the arena is 30 patches wide $(n = 30)$. The different levels of shading represent different densities of hosts and parasitoids: black squares represent empty patches; dark shades becoming paler represent patches with increasing host densities; light shades to white represent patches with hosts and increasing parasitoid densities. (a) spiral: $\mu_H = 1$, $\mu_P = 0.89$; (b) spatial chaos: $\mu_H = 0.2$, $\mu_P = 0.89$; (c) crystalline structures: $\mu_H = 0.05$, $\mu_P = 1$. Case (d) is a similar map obtained with Lotka-Volterra local dynamics, using $\lambda = 1.4$, the parasitoid death rate $d = 0.9$, $\mu_H = 0.8$, $\mu_P = 0.8$; it exhibits highly variable spirals. After [38].

Even richer spatial dynamics can ensue when we extend this work to systems of three or more interacting species. In such systems, coexistence of competing species often appears to be associated with some degree of persistent spatial segregation, even when the environment is uniform. In particular, for two parasitoid species sustained by a single host species, self-organised spatial dynamics can result in one competitively inferior predatory species being confined to small, relatively static "islands" within the habitat. Figure 4 illustrates this phenomenon. The natural ecological interpretation would be that this species is surviving by virtue of isolated pockets of favourable habitat; but, in fact, the substrate is spatially uniform, and the spatial structure seen in Figure 4 is self-organised by nonlinearities in the dynamics [40]

Two conclusions may be drawn from this work. First, local movement in an environment of identical patches can, by itself, enable otherwise unstable communities of interacting populations to persist together. Second, simple and fully deterministic movement rules can give self-organised spatial patterns of may different kinds—spiral waves, spatial chaos, crystal lattices. In particular, the spatial chaos thus generated, within an inherently homogeneous set of patches, can appear indistinguishable from random environmental heterogeneity. For a more complete discussion, see [37, 40, 41].

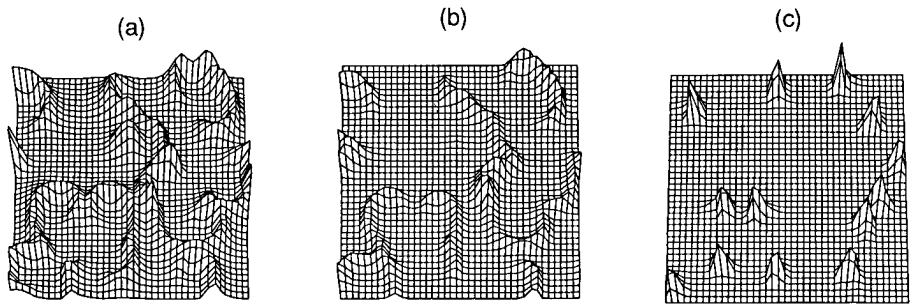

FIGURE 4. Maps of the spatial distribution of density (with linear vertical scales) for: (a) hosts, (b) highly dispersive parasitoids, and (c) relatively sedentary parasitoids. The maps represent a snapshot from the dynamics of a persistent host-parasitoid-parasitoid system of the kind described in the text, with $\lambda = 2$, $\mu_H = 0.5$, $\mu_{P1} = 0.5$, $\mu_{P2} = 0.05$, $a = 1.3$. The grids must be mentally superimposed to perceive the relationships between the densities of the various species. Spiral foci exist at the ends of the "mountain ranges" in (a), excluding the ends at the edges of the grid. In the time-evolution of the system, the mountain ranges are the peaks of population density waves and are thus in continual motion. The foci, by contrast, remain in almost exactly the same place for very long times. The vertical scale in (c) is extended tenfold. After [40].

## SPATIAL CHAOS AND EVOLUTIONARY GAMES

From Darwin's time to ours, one of the central problems of evolution has been the origin and maintenance of cooperative or altruistic behaviour.

Much attention has been given to the Prisoner's Dilemma (henceforth PD), as a metaphor for the problems surrounding the evolution of cooperative behaviour [42]. In its standard form, the PD is a game played by two players, each of whom may choose either to cooperate, $C$, or defect, $D$, in any one encounter. If both players choose $C$, both get a payoff of magnitude $R$; if one defects while the other cooperates, $D$ gets the game's biggest payoff, $T$, while $C$ gets the smallest, $S$; if both defect, both get $P$. With $T > R > P > S$, the paradox is evident. In any one round, the strategy $D$ is unbeatable (being better than $C$ whether the opponent chooses $C$ or $D$). But by playing $D$ in a sequence of encounters, both players end up scoring less than they would by cooperating (because $R > P$). Following Axelrod and Hamilton's pioneering work [43], many authors have sought to understand which strategies (such as "tit-for-tat") do best when the game is played many times between players who remember past encounters. These theoretical analyses, computer tournaments, and laboratory experiments continue, with the answers depending on the extent to which future payoffs are discounted, on the ensemble of strategies present in the group of players, on the degree to which strategies are deterministic or error-prone (e.g., imperfect memories of opponents or of past events), and so on [42, 44–46].

Nowak and May [47–50] have recently given a new twist to this discussion, by considering what happens when the game is played with close neighbours in some 2-dimensional spatial array: "spatial dilemmas".

We consider only two kinds of players: those who always cooperate, $C$, and those who always defect, $D$. No attention is given to past or likely future encounters, so no memory is required and no complicated strategies arise. These memoryless "players"—who may be individuals or organised groups—are placed on a 2-dimensional, $n \times n$ square lattice of "patches"; each lattice-site is thus occupied either by a $C$ or a $D$. In each round of our game (or at each time step or each generation), each patch-owner plays the game with its immediate neighbours. The score for each player is the sum of the payoffs in these encounters with neighbours. At the start of the next generation, each lattice-site is occupied by the player with the highest score among the previous owner and the immediate neighbours. The rules of this simple game among $n^2$ players on the $n \times n$ lattice are thus completely deterministic. Specifically (but preserving the essentials of the PD), we choose the payoffs of the PD matrix to have the values $R = 1$, $T = b (b > 1)$, $S = P = 0$. That is, mutual cooperators each score 1, mutual defectors 0, and $D$ scores $b$ (which exceeds unity) against $C$ (who scores 0 in such an encounter). The parameter $b$, which characterises the advantage of defectors against cooperators, is thus the only parameter in the model; none of the findings are qualitatively altered if we instead set $P = \varepsilon$, with $\varepsilon$ positive but significantly below unity (so that $T > R > P > S$ is strictly satisfied). In the figures that follow, we assume the boundaries of the $n \times n$ matrix are fixed, so that players at the boundaries simply have fewer neighbours; the qualitative character of our results is unchanged if we instead choose periodic boundary conditions (so that the lattice is really a torus). The figures are for the case when the game is played with the eight neighbouring sites (the cells corresponding to the chess king's move) and with one's own site (which is reasonable if the players are thought of as organised groups occupying territory). As amplified below, the essential conclusions remain true if players interact only with the four orthogonal neighbourhoods in square lat-tices or with six neighbours in hexagonal lattices. The results also hold whether or not self-interactions are included [48, 51].

Using an efficient computer program in which each lattice-site is represented as a pixel of the screen, we have explored the asymptotic behaviour of the above-described system for various values of $b$ and with various initial proportions of $C$ and $D$ arranged randomly or regularly on an $n \times n$ lattice ($n = 20$ and up). The dynamical behaviour of the system depends on the parameter $b$; the discrete nature of the possible payoff total means that there will be a series of discrete transition-values of $b$ that lead from one dynamical regime to another. These transition-values and the corresponding patterns are described in detail elsewhere [48]. The essentials, however, can be summarised in broad terms. If $b > 1.8$, a $2 \times 2$ or larger cluster of $D$ will continue to grow at the corners (although not necessarily along the edges, for large squares); for $b < 1.8$, big $D$-clusters shrink. Conversely, if $b < 2$, a $2 \times 2$ or larger cluster of $C$ will continue to grow; for $b > 2$, $C$-clusters do not grow. A particularly interesting regime is therefore $2 > b > 1.8$, where $C$-clusters can grow in regions in $D$ and also $D$-clusters can grow in regions of $C$. As intuition might suggest, in this interesting regime we find chaotically varying spatial arrays, in which $C$ and $D$ both persist in shifting patterns. Although the detailed patterns change from generation to generation—as both $C$-clusters and $D$-clusters expand, collide, and fragment—the asymptotic overall fraction of sites occupied by $C$, $f_C$, fluctuates around 0.318 for almost all starting proportions and configurations.

Color Figure (a) (following page 306) illustrates a typical asymptotic pattern in this regime $2 > b > 1.8$, and shows the typical patterns of dynamic chaos found for almost all starting conditions in this regime. Figure 5a on the next page adds a temporal dimension to Color Figure (a), showing the proportion of sites occupied by $C$ in successive time-steps (starting with 40% $D$). The asymptotic fraction, $f_C$, shown in Figure 5a is found for essentially all starting proportions and configurations, for these $b$-values.

Color Figure (b) is perhaps more in the realm of aesthetics than biology. Again $2 > b > 1.8$, but now we begin ($t = 0$) with a single $D$ at the center of a $99 \times 99$ lattice of $Cs$. The figure shows the consequent symmetrical pattern 200 time-steps later. Such patterns, each of which can be characterised in fractal terms, continue to change from step to step, unfolding a remarkable sequence—dynamic fractals. The patterns show every lace doily, rose window, or Persian carpet you can imagine.

As Figure 5b shows, the asymptotic fraction of $C$ is as for the chaotic pattern typified by Color Figure (a) and Figure 5a. Many of the dynamic features of the symmetric patterns typified by Color Figure (b) can be understood analytically. In particular, we can make a crude estimate of the asymptotic $C$-fraction, $f_C$, for very large such symmetrical patterns by referring to the geometry of the $D$-structure. The $D$-structures are closed-boundary squares in generations that are powers of 2, $t = 2^n$; hence $f_C(t)$ has minima at generations that are powers of 2. These squares now expand at the corners and erode along the sides, returning to square shape after another doubling of total generations. On this basis, a crude approximation suggests that, $i$ time steps en route from $t$ to $2t$, there will be roughly $4(2i)(2t+1-2i)$ $C$-sites within the $D$-structure of size $(2t+1+2i)^2$;

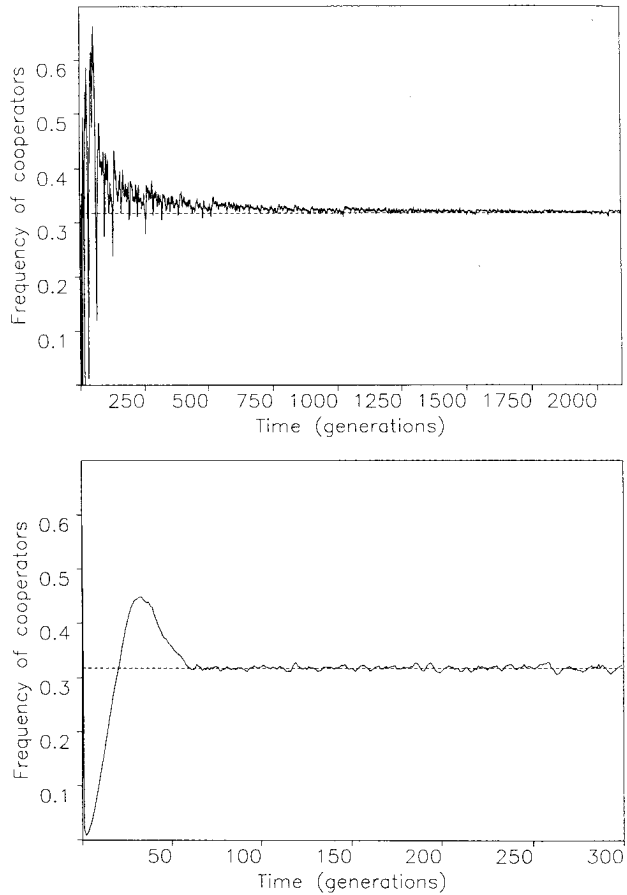

FIGURE 5. (a) The frequency of cooperators, $f_C(t)$, for 300 generations, starting with a random initial configuration of $f_C(0) = 0.6$. The simulation is performed on a $400 \times 400$ square lattice with fixed boundary conditions, and each player interacts with 9 neighbours (including self). The dashed line represents $f_C = 12 \ln 2 - 8 \approx 0.318$ (see text). (b) The frequency of $C$ within the dynamic fractal generated by a single $D$ invading an infinite array of $C$. At generation $t$, the width for the growing $D$-structure is $2t + 1$, and Figure 4b shows the frequency of $C$, $f_C(t)$, within the square of size $(2t+1)^2$ centred on the initial $D$-site, as a function of $t$. Again, the dashed line represents the approximation discussed in the text. After [47, 48].

hence the asymptotic $C$-fraction, $f_C$, for very large such symmetrical patterns is

$$(8) \qquad f_C \approx 4 \int_0^1 s(1-s)(1+s)^{-2} ds = 12 \ln 2 - 8.$$

This approximation, $f_C \approx 0.318 \ldots$, is indicated by the dashed line in Figures 5a and 5b. It agrees with the numerical results for the symmetric case, Figure 5b, significantly better than we would have expected. Why this approximation also works well for the irregular, spatially chaotic patterns, Figure 5a, we do not know [48].

The above work is all based on symmetric lattices in 2-dimensions, deterministic winning (the local winner is always the cell with the highest local score), and discrete time (with all cells updated in unison). In the real world, spatial distributions of cells are likely to be irregular, and winning is likely to have chance elements. There will also be circumstances where interactions occur in continuous, rather than discrete, time, although this is certainly not always the case in problems of biological interest (to the contrary of the extreme view expressed by Hubermann and Glance [52]; see the discussion in [50]). We have therefore generalised studies of "spatial dilemmas" to embrace these complications: spatially random distributions of cells, with the game now being played with neighbours within some specified distance; "probabilistic winning", whereby any relevant neighbour may win a site, with probabilities that depend to a specified extent on the relative scores; and continuous rather than synchronised updating of site ownerships [49, 50].

The essential conclusions of our earlier, and more restricted, analysis remain intact under these generalisations. For a broad band of values in the cheating-advantage parameter $b$, $C$ and $D$ can persist together, in enduring polymorphism. On the other hand, the beautiful "Persian Carpets" of Color Figure (b) of course cannot arise if there is any departure from strict symmetry and/or determinism, as there will be if we have probabilistic winning, or spatial irregularity, or continuous updating, much less an asymmetric starting configuration. As emphasised elsewhere, these entrancing patterns are of aesthetic and mathematical interest but have no direct biological significance (except insofar as mathematical understanding of such special cases can help illuminate the dynamics of the game more generally). For a more detailed discussion, see [48–51].

In short, the PD is an interesting metaphor for the fundamental biological problem of how cooperative behaviour may evolve and be maintained. Essentially all previous studies of the PD are confined to individuals or organised groups who can remember past encounters, who have high probabilities of future encounters (with little discounting of future payoffs), and who use these facts to underpin more-or-less elaborate strategies of cooperation or defection. The range of real organisms obeying these constraints is very limited. In contrast, the spatially embedded models involve no memory and no strategies: the players are pure $C$ or pure $D$. Deterministic interactions with immediate neighbours in 2-dimensional spatial arrays, with success (site, territory) going in each generation to the local winner, is sufficient to generate astonishingly complex and spatially chaotic patterns in which cooperation and defection persist indefinitely. The details of the patterns depend on the magnitude of the advantage accruing to defectors (the value of $b$), but a range of values leads to polymorphic patterns, whose nature is almost always independent of the initial proportions of $C$ and $D$.

These studies suggest that deterministically generated spatial structure within populations may often be crucial for the evolution of cooperation, whether it be among molecules, cells, or organisms. Other evolutionary games, such as Hawk-Dove [53], which recognise such chaotic or patterned spatial structure may be more robust and widely applicable than those that do not. More generally, such self-generated and complex spatial structures may be relevant to the dynamics of a wide variety of spatially extended systems: Turing models, 2-state Ising models, and various models for prebiotic evolution (where it seems increasingly likely that chemical reactions took place on surfaces rather than in solutions).

## Concluding remarks

In Stoppard's [54] latest play, *Arcadia*, one of the characters—a graduate student working on the irregular cycles observed in the abundance of grouse—offers an extreme view of the significance of chaos: "The unpredictable and the predetermined unfold together to make everything the way it is. It's how nature creates itself, on every scale, the snowflake and the snowstorm. It makes me so happy . . . . A door like this has cracked open five or six times since we got up on our hind legs. It's the best possible time to be alive, when almost everything you thought you knew is wrong."

While such a statement may represent dramatic licence, I think it true that recent advances in understanding chaotic dynamics (and associated fractal geometries) really are opening new doors. And not only for researchers in mathematics and the biological and physical sciences. A lot of the basic mathematics has the combination of simplicity and surprise that enables an undergraduate, or even a high school student, to be taken to the frontiers of this new subject and to enjoy what is found there. The sooner this relevant material is incorporated in undergraduate and high school curricula, the better for all of us.

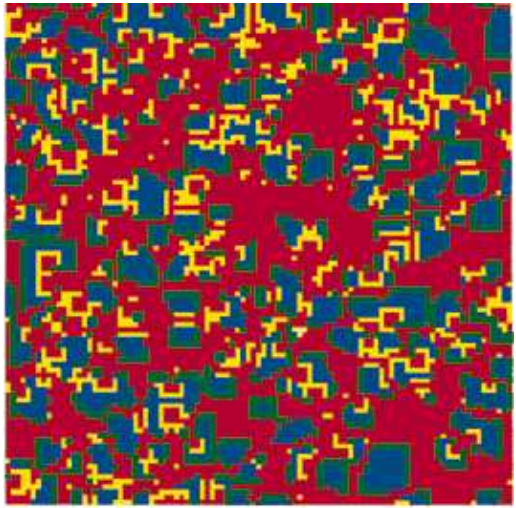

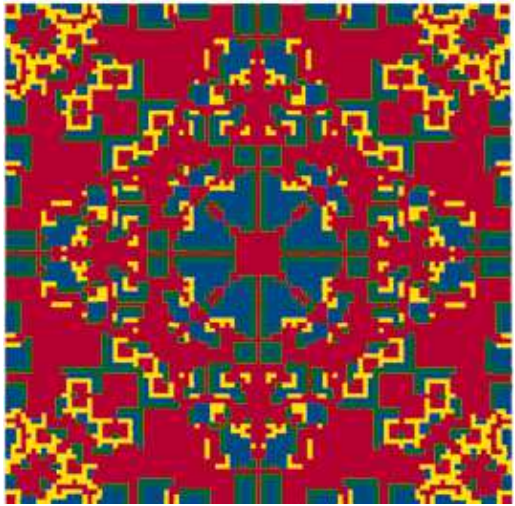

Figure The spatial Prisoner's Dilemma can generate a large variety of qualitatively different patterns, depending on the magnitude of the parameter, $b$, which represents the advantage for defectors. Both (a) and (b) are for the interesting region when $2 > b > 1.8$. The pictures are coded as follows: blue squares represent a cooperator ($C$) that was already a $C$ in the preceding generation; red (the commonest square) represents a defector ($D$) following a $D$; green is a $C$ following a $D$; and yellow shows a $D$ following a $C$. (a) This simulation is for a $99 \times 99$ lattice with fixed boundary conditions, starting with a random configuration with 10% $D$ and 90% $C$. The figure shows the asymptotic pattern after 200 generations: spatial chaos. (b) Beautiful "fractal kaleidoscopes" ensue if the initial pattern is symmetric (the rules preserve such symmetry). Here the simulation is started with a single $D$ at the centre of a $99 \times 99$ field of $C$ with fixed boundary conditions. The figure shows that pattern 200 generations later. After [47, 48].

Department of Zoology, University of Oxford, Oxford OX1 3PS, UK, and Imperial College, London
*E-mail address*: `rmay@vax.oxford.ac.uk`